\documentclass[12pt,a4paper]{amsart}

\newtheorem{theo+}              {Theorem}           [section]
\newtheorem{prop+}  [theo+]     {Proposition}
\newtheorem{coro+}  [theo+]     {Corollary}
\newtheorem{lemm+}  [theo+]     {Lemma}
\newtheorem{exam+}  [theo+]     {Example}
\newtheorem{rema+}  [theo+]     {Remark}
\newtheorem{defi+}  [theo+]     {Definition}

\newenvironment{theorem}{\begin{theo+}}{\end{theo+}}
\newenvironment{proposition}{\begin{prop+}}{\end{prop+}}
\newenvironment{corollary}{\begin{coro+}}{\end{coro+}}

\usepackage{amsthm}
\theoremstyle{plain} \theoremstyle{remark}

\newtheorem{example}{Example}

\newtheorem*{ack}{\bf Acknowledgments}

\def\E{/\kern-1.0em \equiv }

\author{}

\begin{document}
\title[Biharmonic  conformal maps]{Some recent work on biharmonic  conformal maps}
\author{Ye-Lin Ou $^{*}$}
\thanks{$^{*}$ This work was supported by a grant from the Simons Foundation ($\#427231$, Ye-Lin Ou).}
\address{Department of
Mathematics,\newline\indent Texas A $\&$ M University-Commerce,
\newline\indent Commerce, TX 75429,\newline\indent USA.\newline\indent
E-mail:yelin.ou@tamuc.edu }
\date{09/10/2019}
\maketitle
\section*{Abstract}
\begin{quote} This note  reviews  some of the  recent work on biharmonic conformal maps (see \cite{OC}, Chapter 11, for a detailed survey). It will be focused on biharmonic conformal immersions and biharmonic conformal maps between manifolds of the same dimension and their links to isoparametric functions and Yamabe type equations, though biharmonic morphisms (maps that preserve solutions of bi-Laplace equations),  generalized harmonic morphisms (maps that pull back germs of harmonic functions to germs of biharmonic functions), and biharmonic conformal and Riemannian submersions will also be touched.
{\footnotesize } 
\end{quote}

\section{Biharmonic conformal maps}
\indent A {\it weakly conformal map}  is a conformal immersion $\phi: (M^m, g)\to(N^n,h)$ away from its branch points where  ${\rm d} \phi_p=0$. So, for a weakly conformal map, we have $\phi^*h=\lambda^2g$ for a function  $\lambda\ge 0$ called the conformal factor of $\phi$, and  $\lambda(p)=0$ if and only if $p\in M$ is a branch point of $\phi$. A weakly conformal map with conformal factor $\lambda>0$ (respectively, $\lambda=C>0$, or $\lambda=1$) is called a {\it conformal immersion} (respectively, a {\it homothetic}, or an {\it isometric immersion}).\\
\indent A {\em horizontally weakly conformal map} is a map $\phi:(M,
g)\to (N, h)$  such  that for any $p \in M$, either  ${\rm d} \phi_p=0$, or
 ${\rm d} \phi_{p}|_{ \mathcal{H}_p} : \mathcal{H}_p \to T_{\phi(p)} N$ is
conformal and onto, where $ \mathcal{H}_p=({\rm ker} {\rm d} \phi_p)^{\bot}$ is the horizontal subspace of $\phi$.  It is well known (see, e.g., \cite{Baird03}) that for a horizontally weakly conformal map, there exists a function $\lambda\ge 0$, called the {\it dilation} of $\phi$, such that $h({\rm d} \phi
(X),{\rm d} \phi (Y)) = \lambda^{2}(p)g(X,Y)$ for any horizontal vectors $X, Y \in
 \mathcal{H}_p$.  Note that  $\lambda(p) = 0$ if and only if ${\rm d} \phi_{p} = 0$. A horizontally weakly conformal map with dilation $\lambda>0$ (respectively, $\lambda=C>0$, or $\lambda=1$) is called a {\it conformal submersion} (respectively, a {\it homothetic}, or a {\it Riemannian submersion}).\\
\indent Clearly, we have the following relationships among various conformal maps:
\begin{center} {\scriptsize $
\{ {\rm Isometric\; immersions}\} \subset  \{ {\rm Conformal\; immersions}\} \subset  \{ {\rm Weakly\; conformal\; maps}\}$;\\
$\{ {\rm Riemannian\; submersions}\}\subset  \{ {\rm Conformal\; submersions}\}\subset  \{ {\rm Horizontally\; weakly\; conformal\; maps}\}.
$}\end{center}

\indent When the domain and the target manifolds have the same dimension, the weak conformality agrees with the horizontal weak conformality. In this case,   conformal immersions are precisely  conformal submersions, and we call them {\it conformal maps} in this note.  Note that, unlike most other books,  our definition of conformal maps does not require conformal maps be diffeomorphisms, although they are  local diffeomorphisms by the well known inverse mapping theorem.\\
\indent A {\it harmonic map} is a map $\phi: (M^m, g)\to(N^n,h)$ between Riemannian manifolds whose tension field $\tau (\phi)={\rm Trace}_g\nabla{\rm d}\phi$ vanishes identically. Harmonic maps include harmonic functions, geodesics, minimal submanifolds, and Riemannian submersions with minimal fibers as special cases.  We refer the reader to \cite{EL78}, \cite{EL83}, \cite{EL88} and \cite{YS97} for the beautiful theory, important applications and interesting links of harmonic maps and the related topics.\\
\indent A {\it biharmonic map} is a map  $\phi: (M^m, g)\to(N^n,h)$ between Riemannian manifolds that is a critical point of the bienergy functional
$E_2(\phi, g)=\frac{1}{2}\int_{M} |\tau (\phi)|^2dv_g\,,$
where $\tau (\phi )$ is the tension field of $\phi$. Biharmonic map equation is a system of the 4th order elliptic PDEs (see \cite{Ji86}):  
\begin{equation} \label{bihar}
\tau_2(\phi ) :={\rm Trace}_g\left(  (\nabla^{\phi})^2 \tau (\phi ) - R^N( {\rm d}\phi, \tau (\phi )){\rm d}\phi\right) = 0
\end{equation}
where $R^N$ denotes the Riemannian curvature of the target manifold $N$.\\
\indent A {\it biharmonic submanifold} is one whose defining isometric immersion is a biharmonic map. From  the definition of biharmonic maps and the well known fact that a submanifold is minimal if and only if its defining isometric immersion is harmonic, we have the following relations:
 {\scriptsize \begin{eqnarray}\notag
\{ {\rm Harmonic\; maps}\} &\subset & \{ {\rm Biharmonic\; maps}\};\\\notag
\{ {\rm Minimal\; submanifolds}\}&\subset & \{ {\rm Biharmonic\; submanifolds}\}.
\end{eqnarray}}
 For this reason, we use  the name {\it proper biharmonic maps} for those  biharmonic maps which are  not harmonic, and {\it proper biharmonic submanifolds} for those  biharmonic submanifolds which are  not minimal.\\
 \indent By the nature of the biharmonic map equation in general, it is very difficult to get  a solution, or try to get some geometric and/or topological information of the spaces which support such maps. However, some recent study have shown that studying biharmonic maps under some geometric constraints is both practical and fruitful. For example, by adding geometric condition of ``being isometric immersion", we have
 {\scriptsize 
 \begin{eqnarray}\notag
 \{ {\rm Biharmonic\; submanifolds}\}= \{ {\rm Biharmonic\; maps}\}\cap \{ {\rm Isometric\; immersions}\},
 \end{eqnarray}}
 which have been receiving a growing interest from many researchers and some interesting progress has been made  since 2000. We refer the reader to  a recent survey \cite{Ou16} for a progress report on the study of biharmonic submanifold.\\
 \indent Another natural  geometric condition to  impose is ``being  horizontally weakly conformal (HWC for short)", thus we study biharmonic HWC maps, which contains the following interesting subsets 
 {\scriptsize  \begin{eqnarray}\notag  
 \{ {\rm Bihar.\; Riemannian\; submersions}\}\subset \{ {\rm Bihar.\; Conformal\; submersions}\}\subset \{ {\rm Bihar.\; WHC\; maps}\}.
 \end{eqnarray}}
 Motivations to study biharmonic HWC maps include:
 \begin{itemize}
 \item [(i)] Bihar. morphisms  are special cases of bihar. HWC maps;
 \item[(ii)] Generalized har. morphisms into $\mathbb{R}^n$  are special bihar. HWC maps;
  \item[(iii)] Bihar. HWC maps are generalization of har. morphisms since
  {\scriptsize  \begin{eqnarray}\notag  
 \{ {\rm Har.\; morphisms}\}= \big( \{ {\rm Har.\; maps}\}\cap \{ {\rm  WHC\; maps}\}\big)\subset \big( \{ {\rm Bihar.\; maps}\}\cap \{ {\rm WHC\; maps}\}\big).
 \end{eqnarray}}
   \end{itemize}
\indent Recall that  a harmonic morphism is a map between Riemannian manifolds that preserves solutions of Laplace equation. Harmonic morphisms are characterized as HWC harmonic maps (\cite{Fu78}, \cite{Is79}). We refer the reader to the book \cite{Baird03} for a comprehensive study of the geometry of harmonic morphisms.\\ 
\indent Recall also that a {\it biharmonic morphism} (\cite{Ou00}) is a map  $\phi: (M, g)\to (N, h)$ between Riemannian manifolds which  preserves solutions of bi-Laplace equation in the sense that  $\Delta^2_N\,f=0$ implies $\Delta^2_M(f\circ \phi)=0$ for any function $f:N\supset V\to \mathbb{R}$ with $\phi^{-1}(V)$ non-empty. It was proved in \cite {Ou00}, \cite{LO01} and \cite{LO10} that any biharmonic morphism is a biharmonic HWC map.\\
\indent A  more general type of morphisms of harmonicity  called {\it generalized harmonic morphisms} were recently studied in \cite{GO18}. By  definition, a generalized harmonic morphism is a map between Riemannian manifolds that pulls back germs of harmonic functions to germs of biharmonic functions. Among other things, it was proved in \cite{GO18} that
 \begin{itemize}
 \item [(i)] $ \{ {\rm Har.\; morphisms}\} \subset  \{ {\rm Generalized\; har.\; morphisms}\}$
 \item[(ii)] $ \{ {\rm Bihar.\; morphisms}\} \subset  \{ {\rm Generalized\; har.\; morphisms}\}$
\item[(iii)] Generalized har. morphisms into $\mathbb{R}^n$  are special bihar. HWC maps;
\item[(iv)] There are many examples of proper biharmonic Riemannian submersions.
 \end{itemize}
   \vskip0.1cm 
{\it Biharmonic Riemannian submersions} were first studied in \cite{Oni02}. For some recent work on this subject, see  \cite{BMO07}, \cite{LO10}, \cite{WO11}, \cite{GO18}, \cite{Ur18}, and \cite{AO19}.

\section{Biharmonic conformal maps between manifolds of the same dimension}

{\bf 2.1 Biharmonic conformal map equations}\\
 It is easily checked that the tension field  of a conformal map  $\phi: (M^m, g)\to (N^m, h)$ is $\tau (\phi ) = - (m-2) d\phi (\nabla \ln \lambda )$. Thus, the picture of harmonic maps between manifolds of the same dimension is very simple: every conformal map $\phi: (M^2, g)\to (N^2, h)$ is harmonic;  and for $m\ne 2$, a conformal map is $\phi: (M^m, g)\to (N^m, h)$ harmonic  if and only if  $\lambda$ is a constant, i.e., $\phi$ is a homothety.\\
 \indent The following theorem shows that biharmonic conformal maps are far more complicated and hence have much more to study.
\begin{theorem}\cite{BO2}\label{MT3}
A conformal map $\phi :(M^m, g) \to (N^m, h)$ with  $m\ge 3$ is biharmonic if and only if the conformal factor $\lambda$ solves the PDE:
\begin{eqnarray}\label{SnD}
\lambda\,\nabla\,\Delta \lambda -3(\Delta \lambda)\nabla \lambda  -\frac{m-4}{2}\nabla |\nabla \lambda|^2+2\lambda {\rm Ric}^M(\nabla \lambda)=0,
\end{eqnarray}
where and in the rest of the paper, $\nabla$, $|, |$, and
$\Delta$ denote the gradient, the norm, and the Laplacian defined by the metric $g$, 
In particular,  a conformal map  from an Einstein $4$-manifold with ${\rm Ricci}^M = a g$ is  biharmonic if and only if
\begin{equation}\label{4D}
\Delta \lambda - a\lambda = A \lambda^3,
\end{equation}
where $A$ is a constant such taht
\begin{equation} \label{A-cf-scal}
6A+\frac{2a}{\lambda^2} + R_h = 0
\end{equation}
holds for the scalar curvature $R_h$ of the target manifold.
\end{theorem}

{\bf 2.2 Interesting links to isoparametric functions and Yamabe type equations}\\
\noindent{\bf (i)  Biharmonicity  and isoparametric functions.} \\The following interesting link was found in \cite{Baird08}
\begin{theorem} {\rm (\cite{BK}, \cite{Baird08}, see also \cite{BO2})}
 For $m>2$ and $m \neq 4$, a conformal map $\phi :(M^m, g) \to (N^m, h)$ from an Einstein manifold is biharmonic if and only if the conformal factor $\lambda$ is an isoparametric function.
\end{theorem}
\indent Here, an isoparametric function is a function  $f: M\to R$ on a Riemannian manifold  such that  ${\rm (i)}\; |\nabla\, f|^2=\alpha (f)$, and ${\rm (ii)} \;\Delta f=\beta(f)$ for some real variable functions $\alpha$ and $\beta$. Isoparametric functions on space forms have rich geometry related to hypersurfaces of constant principal curvatures  which have been studied since 1918. For a good survey on this subject, see \cite{Th00}.\\
\noindent {\bf (ii) Biharmonic equation  is a Yamabe type Equation.}  It is well known that the Yamabe problem is to find a metric  $h$, within the conformal class $ [g] =\{ \sigma g| \, \sigma: M\to (0, \infty)\}$ of metrics on a closed Riemannian manifold $(M^m, g)$, that  has constant scalar curvature $R_h$. For $m \geq 3$ and  $h=\lambda^{\frac{4}{m-2}}g$, the problem is equivalent to solving  the Yamabe equation: 
\begin{equation}\label{Ya}
-\frac{4(m-1)}{m-2}\Delta_g\,\lambda+R_g\,\lambda=R_h\,\lambda^{\frac{m+2}{m-2}},\;\;R_h={\rm constant}.
\end{equation}
In Particular, when $m=4$, the Yamabe equation reads:
\begin{equation}\label{Ya}
\Delta_g\,\lambda-\frac{1}{6}R_g\,\lambda=-\frac{1}{6}R_h\,\lambda^3,\;\;R_h={\rm constant}.
\end{equation}

More generally,  the following {\it Yamabe type equation}  
\begin{equation}\label{gYa}
\Delta_g\,u-\alpha\,u=-\beta\,u^3, \;\;\alpha\; {\rm and}\;\; \beta\;\;{\rm are \;functions},
\end{equation}
has been studied by many mathematicians (see, e.g., \cite{Au76a}, \cite{Au76b}, \cite{Au98}, and \cite{VW17}). \\
\indent Theorem \ref{MT3} shows that biharmonic equation for conformal maps from an Einstein $4$-manifold  is a Yamabe type equation.\\
\noindent {\bf (iii) Biharmonicity and Sobolev embedding Theorem.} \\ The following link was found in \cite{BO2}
\begin{corollary}
For $m=4$, the Euler-Lagrange equation for the extremal functions saturating the Sobolev inequality
\begin{equation}\notag
c\left(\int_{\mathbb{R}^m}|v|^{2m/(m-2)}dx   \right)^{(m-2)/m}\le \int_{\mathbb{R}^m}|\nabla v|^2dx
\end{equation}
 is 
\begin{equation}\label{EL4}
\Delta v=-2v^3,\hskip1cm {\rm on}\;\; \mathbb{R}^4,
\end{equation}
which is a special case of the biharmonic equation for conformal maps from Euclidean space $\phi:\mathbb{R}^4\supseteq U\to (N^4, h)$. \end{corollary}
We refer the reader to \cite{Ch} for an explanation of the Euler-Lagrange equation for the extremal functions saturating the Sobolev inequality.\\

{\bf 2.3 Some recent results}\\
\indent We would like to point out that  it is the discovery of the  link  between biharmonic equation for conformal maps  and the Yamabe type equation and he Euler-Lagrange equation for the extremal functions saturating the Sobolev inequality that provides us  very useful tools in our classification of biharmonic M\"obius transformations, understanding of biharmonic metrics, and proving some existence theorems  about biharmonic conformal maps from a compact Einstein $4$-manifold. The following is a  list of some results from \cite{BO2}, we refer the reader to the original paper for all details and other results.

$\clubsuit$  A  conformal map $\phi :(M^4, g) \to (N^4, h)$ with $\phi^{*}h=\lambda^2g$ from an Einstein manifold with ${\rm Ricci}^M = ag$ into a manifold with constant scalar curvature is biharmonic if and only if it is a homothety (hence a harmonic map), or
\begin{equation}\notag
a=0,\;\; {\rm and}\;\;\Delta \lambda = -\frac{R_h}{6} \lambda^3.
\end{equation}
\indent $\clubsuit$  The only possible  case for proper biharmonic conformal map to exist (even locally) between $4$-dimensional  space forms is: $\mathbb{R}^4\supseteq U\to M^4(C)$, where $M^4(C)$ denotes a $4$-dimensional space form. This confirms and extends a result in \cite{MOR15} which states that the only rotationally symmetric, proper biharmonic conformal diffeomorphisms between $4$-dimensional space form appear in the case of $\mathbb{R}^4\supseteq U\to M^4(C)$. \\
\indent $\clubsuit$ A conformal map $\phi: \mathbb{R}^4\supseteq U\to \mathbb{R}^4$, or  $\phi: \mathbb{R}^4\supseteq U\to S^4$  is biharmonic if and only if it is a restriction of a M\"obius transformation $\phi: \mathbb{R}^4\supseteq U\to \mathbb{R}^4$ with
\begin{equation} \notag
\phi(x)= a + \frac{\alpha A(x-b)}{|x-b|^{\varepsilon}} \quad a, b \in \mathbb{R}^4, \ \alpha \in \mathbb{R} \setminus\{ 0\} , \ A \in O(4),\ \varepsilon \in \{ 0, 2\},
\end{equation}
for all $x \in \mathbb{R}^4$. The first map is proper biharmonic if and only if  $\varepsilon = 2$ whilst the seond map is always proper biharmonic. \\
\indent $\clubsuit$ There exists a function $\lambda$ on a closed Einstein $4$-manifold with negative Ricci curvature so that  $1_g: (M^4, g)\to (M^4, h=\lambda^2g)$ is a proper biharmonic map.\\
\indent $\clubsuit$  There exists no proper biharmonic conformal map from a close Ricci flat $4$-manifold  into any Riemannian $4$-manifold;\\
\indent $\clubsuit$ There exists an infinite family of conformal metrics $\lambda^2g$ on the Euclidean $4$-sphere  so that the identity map $1_g: (S^4, g)\to (S^4, h=\lambda^2g)$ is a proper biharmonic conformal diffeomorphism. Note that by the 2nd item of this list, no (not even local) proper biharmonic conformal map exists between the Euclidean $4$-spheres.

Recall that a Riemannian metric ${\bar g}$ on a Riemannian manifold $(M^m, g)$ is called a {\it harmonic metric} (\cite{CN84}) (respectively, {\it biharmonic metric}) with respect to $g$ if the identity map $1_g: (M^m, g)\to (M^m, {\bar g})$ is a harmonic map (respectively, {\it biharmonic map}). 

\begin{corollary}{\rm (\cite{Ba04}, \cite{BO2})}
 The identity map $1_g: (M^m, g)\to (M^m, {\bar g}=\lambda^{2}g)$ is  biharmonic (i.e., the conformally related metric ${\bar g}=\lambda^{2}g$ is biharmonic  with respect to $g$)  if and only if 
\begin{eqnarray}\label{SFI}
\lambda\,\nabla\,\Delta \lambda -3(\Delta \lambda)\nabla \lambda  -\frac{m-4}{2}\nabla |\nabla \lambda|^2+2\lambda {\rm Ric}^M(\nabla \lambda)=0.
\end{eqnarray}
\end{corollary}
 From this we have:

$\clubsuit$  A  metric $h=\lambda^2 dx^2$ on an open set $U\subseteq \mathbb{R}^4$ is biharmonic if and only if $\lambda$ is a solution of the Yamabe equation or, equivalently,  $h=\lambda^2 dx^2$ has constant scalar curvature $R_h$ and $\Delta \lambda = -\frac{R_h}{6} \lambda^3$ (\cite{BO2}).

$\clubsuit$ The conformally flat metric $h=|x|^{2\alpha} {\rm d}x^2$ on $\mathbb{R}^4\setminus \{0\}$ is biharmonic with respect to Euclidean metric ${\rm d}x^2$ if and only if $\alpha=-1, \; {\rm or}\;-2$  (\cite{BO2}).

\begin{example}
The  conformally flat metric $\lambda^2 {\rm d} x^2$ with $\lambda= |x|^{-1}, |x|^{-2}, \frac{2}{1\pm |x|^2}$ is biharmonic with respect to Euclidean metric ${\rm d} x^2$ on an open set of $U\subset \mathbb{R}^4$. Equivalently, each of the following identity maps
\begin{eqnarray}\notag
&& 1: (\mathbb{R}^4\setminus \{0\}, {\rm d}x^2) \to \left(\mathbb{R}^4\setminus \{0\},\frac{{\rm d}x^2}{|x|^2}\right),\;\;\; 1: (\mathbb{R}^4\setminus \{0\}, {\rm d}x^2)\to \left(\mathbb{R}^4\setminus \{0\},\frac{{\rm d}x^2}{|x|^4}\right),\\\notag
&&  1: (\mathbb{R}^4, {\rm d}x^2)\to \left(\mathbb{R}^4,\frac{4 {\rm d}x^2}{(1+|x|^2)^2}\right),\hskip1.5cm 1: (\mathbb{B}^4, {\rm d}x^2) \to \left(\mathbb{B}^4, \frac{4 {\rm d}x^2}{(1-|x|^2)^2}\right),
\end{eqnarray} is  a proper biharmonic map (\cite{BO2} and \cite{LO10}).
\end{example}

Note that several results (including the following one) on some obstructions on the existence of biharmonic metric on a compact Einstein manifold have been obtained in \cite{BLO11}.
\begin{theorem}\cite{BLO11}
There exists no proper biharmonic metric on a compact Einstein manifold $(M^m, g)$ with ${\rm Ric}^g=a\,g$ for either (i) $a < 0$ and $m\ge 6$, or (ii) $a\le 0$ and $m>6$.
\end{theorem}

\section{Biharmonic conformal immersions}
\indent In this section, we review some recent work on  biharmonic conformal immersions done in  \cite{Ou09},  \cite{Ou14}, \cite{Ou15}, and \cite{Ou17}.\\
\noindent {\bf 3.1 Biharmonic conformal immersions}\\
\indent It is easily checked (see e.g., \cite{Ta94}) that the tension field of a conformal immersion $\phi : (M^{m},g) \to (N^{n},h)$ with $\phi^*h=\lambda^2g$ is given by  $\tau (\phi ) =m\,\lambda^2\eta - (m-2) {\rm d}\phi (\nabla \ln \lambda )$, where $\eta$ is the mean curvature vector field of the associated  isometric immersion (submanifold) $\phi : (M^{m}, {\bar g}=\lambda^2g) \to (N^{n},h)$. Using this with $m=2$, we have the well-known \\
{\bf Fact: } A  conformal immersion $\phi : (M^2, g) \to (N^n, h)$ of a surface  is harmonic if and only if its image $\phi(M^2)\subset (N^n, h)$ is a minimal surface.\\
\indent It is this interesting link, together with the beautiful theory and important applications  of minimal surfaces, and the relation
{\scriptsize\begin{center}$ \{ {\rm Har.\; conformal\;immersions}\} \subset \{ {\rm Bihar.\; conformal\; immersions}\},$\\\end{center}}
\noindent that motivates us to study biharmonic conformal immersions.

An idea here is to express the biharmonicity of a conformal immersion by using the data (including the Laplacian, the gradient, the shape operator and the mean curvature, etc) of its associated isometric immersion. A reason for this is not only that we can use some of the tools and results from  submanifold theory but also that (as we hope)  the study of  biharmonic conformal immersions will help to understand the geometry of submanifolds.

\indent  The following transformations of the Jacobi operator and the bitension field under a conformal change of the domain metric for a generic map play an important role in  our study of conformal immersions. Note that the idea of using a conformal change of the domain metric of a harmonic map to construct  proper biharmonic maps or metrics had been successfully used in \cite{BK}.

\begin{theorem}\label{confC}\cite{Ou09}
For  a map $\phi : (M^{m},g) \to (N^{n},h)$ and a conformal change of the metric ${\bar g}=\lambda^{2}g$, let $J^{\phi}_{g}$ and $\tau_{2}(\phi,g)$ (respectively,   $J^{\phi}_{\bar g}$ and  $\tau_{2}(\phi,{\bar g})$) denote the Jacobi operator and the bitension field of $\phi$ with respect to $g$ (respectively, ${\bar g}$) respectively. Then, the following transformation hold
\begin{eqnarray}\label{JF}\notag
J^{\phi}_{\bar g}(X)&=& \lambda^{-2}[J^{ \phi}_{g}(X)-(m-2)\nabla^{
\phi}_{\nabla \ln \lambda}X],\;\;{\rm and}\\\notag
\tau_{2}(\phi,{\bar g})&= &\lambda^{-4}\{\tau_{2}(\phi, g)-(m-2)J^{
\phi}_{g}({\rm d}{\phi}(\nabla \ln \lambda))\\\notag && -
2(\Delta {\rm ln}\lambda+(m-4)\left| \nabla \ln \lambda \right|^2)\tau(\phi,
g)+(m-6)\nabla^{ \phi}_{\nabla \ln \lambda}\,\tau(\phi,
g)\\\notag && -2(m-2)(\Delta {\rm ln}\lambda+(m-4)\left|\nabla \ln \lambda \right|^2){\rm d}{\phi}(\nabla \ln \lambda)\\\notag
&&+(m-2)(m-6)\nabla^{ \phi}_{\nabla \ln \lambda}\,{\rm
d}{\phi}(\nabla \ln \lambda)\},
\end{eqnarray}
where the Jacobi operator is defined as
\begin{equation}\label{JO}
J^{\phi}_{g}(X)=-\{{\rm
Trace}_{g}(\nabla^{\phi}\nabla^{\phi}-\nabla^{\phi}_{\nabla^{M}})X
- {\rm Trace}_{g} R^{N}({\rm d}\phi, X){\rm d}\phi\}
\end{equation}
for any vector field $X$ along the map $\phi$ so that $J^{\phi}_{g}(\tau(\phi))=-\tau_2(\phi)$.
\end{theorem}

\begin{proposition}\cite{Ou09}\label{CIso}
A conformal immersion $\phi : (M^{m},g) \to (N^{n},h)$ with the conformal factor $\lambda$ is biharmonic if and
only if
\begin{eqnarray}\notag
\lambda^{4}\tau_{2}(\phi,{\bar g})= -(m-2)J^{ \phi}_{g}({\rm
d}{\phi}(\nabla\,\ln\lambda)) -m(\Delta \lambda^2)\eta+m(m-6)\lambda^2\nabla^{ \phi}_{\nabla\,\ln\,\lambda}\,
\eta,
\end{eqnarray}
where $\tau_{2}(\phi,{\bar g})$ and $\eta$ are the bitension field and the mean curvature vector of the associated submanifold.
\end{proposition}
Using this, we can check the following 
\begin{example}\label{E1}
The map $\phi : (\mathbb{R}^3\times \mathbb{R}^{+}, g=\delta_{ij}) \to
(H^5=\mathbb{R}^4\times \mathbb{R}^{+},h=y_5^{-2}\delta_{\alpha\beta})$ with $\phi(x_1,\ldots,x_4)=(1,x_1,\ldots,x_4)$ gives  a proper biharmonic conformal immersion of the Euclidean space into the Hyperbolic space. 
\end{example}
Note that a similar problem of studying biharmonicity of immersions via conformal changes
on the target manifold was done in \cite{Oni03}, where, among other things,  the author proved that the equator of a Euclidean sphere (viewed as a totally geodesic submanifold) can be turned into a proper biharmonic submanifold via some particular conformal change of the Euclidean metric on the ambient sphere.

{\bf 3.2 Biharmonic conformal immersions of surfaces}\\

For a conformal immersion of a surface we have
\begin{theorem}\label{NEW}\cite{Ou09, Ou15, Ou17}
A conformal immersion $\phi : (M^{2},{\bar g}) \to
(N^3, h)$ with conformal factor $\lambda$ is biharmonic if and only if 
\begin{equation}\label{M03}
\begin{cases}
\Delta (\lambda^2 H )-(\lambda^2H)[|A|^2-{\rm Ric}^N(\xi,\xi)]=0,\\A({\rm grad} (\lambda^2 H))+ (\lambda^2 H) [{\rm grad}
H- \,({\rm Ric}^N\,(\xi))^{\top}]=0.
\end{cases}
\end{equation}
where $\xi$, $A$, and $H$ are the unit normal vector field, the
shape operator, and the mean curvature function  of
the surface $\phi(M)\subset (N^3, h)$ respectively, and the
operators $\Delta,\; {\rm grad}$ and $|,|$ are taken with respect to
the induced metric $g=\phi^{*}h=\lambda^{2}{\bar g}$ on the
surface.
\end{theorem}

It was shown in \cite{Ou10} that an isometric immersion (i.e., hypersurface) $M^m\to (N^{m+1}, h)$ is biharmonic if and only if 
\begin{equation}\label{Bhyper}
\begin{cases}
\Delta H-H |A|^{2}+H{\rm
Ric}^N(\xi,\xi)=0,\\
 2A\,({\rm grad}\,H) +\frac{m}{2} {\rm grad}\, H^2
-2\, H \,({\rm Ric}^N\,(\xi))^{\top}=0.
\end{cases}
\end{equation}

By comparing Equations (\ref{Bhyper}) and (\ref{M03}), we see how biharmonic equation transforms under the change of the conformal factor  from constant $1$ to a function $\lambda$. The similar form of the first equation in the system allows us to use some of the tools used in and results from the biharmonic submanifold theory to study biharmonic conformal immersions.

In general, we have 
{\scriptsize\begin{center}
$\{ {\rm Isometric\; immersions}\}\subset \{ {\rm Conformal\; immersions}\}$ and hence\\
$\{ {\rm Bihar.\; isometric\; immersions}\}=\{ {\rm Bihar.\; submanifolds}\}\subset \{ {\rm Bihar. \;conformal\; immersions}\}$\end{center}}
When the ambient space is a Euclidean space, a well known result,  proved by Chen-Ishikawa \cite{CI91} and Jiang \cite{Ji87} independently, shows that there exists no proper biharmonic isometric immersion from a surface into Euclidean space $\mathbb{R}^3$. However, the following result shows that there are infinitely many proper biharmonic conformal immersions of surfaces into $\mathbb{R}^3$.

\begin{proposition}\label{EZHU}\cite{Ou09}
For constants $C_1$ and  $C_2$ with $C_2/C_1< 0$,  let $\lambda^2=C_1e^{ z}+C_2e^{- z}$
and $D=\{ (\theta,z)\in (0,2\pi)\times \mathbb{R}: z\ne \frac{1}{2}\ln
(-C_2/C_1)\}$. Then, the map $\phi:( D, g=\lambda^{-2}(
d\theta^2+dz^2))\to
\mathbb{R}^3$ with $\phi(\theta,z)=(\cos\,\theta, \sin\,\theta, z)$ is a  proper biharmonic conformal
immersions of  a circular cylinder of radius $1$ into Euclidean space $\mathbb{R}^3$. In particular, for $D=\{ (\theta,z)\in (0,2\pi)\times \mathbb{R}: z\ne 0\}$, we have a proper biharmonic conformal immersion 
$\phi:( D, g=\sinh^{-1}(z)(
d\theta^2+dz^2))\to
\mathbb{R}^3$ with $\phi(\theta,z)=(\cos\,\theta, \sin\,\theta, z)$.
\end{proposition}

The following theorem gives some obstructions to the existence of proper biharmonic conformal immersions.
\begin{theorem}\label{MT2}\cite{Ou17}
Let $\phi : (M^{2},{\bar g}) \to (\mathbb{R}^3, h_0)$ be a
 biharmonic conformal immersion into $3$-dimensional Euclidean space with $\phi^{*}h_0=\lambda^2{\bar g}$. Suppose the associated surface $\phi : (M^{2}, g=\lambda^2{\bar g}) \to (\mathbb{R}^3, h_0)$ is complete with the mean curvature function $H$ satisfying $\int_M\,\lambda^4H^2 \,dv_{g}<\infty$.  Then the biharmonic conformal immersion $\phi$ is a minimal (i.e., harmonic) immersion.
\end{theorem}

A practical  (also relatively easier) question to ask  is: what submanifold admits a biharmonic conformal immersion into its ambient space? Here, a submanifold of $(N^{n}, h)$ given by an isometric immersion $\phi: (M^{m}, g) \to (N^{n},h)$ is said to {\em admit a biharmonic conformal immersion into its ambient space}, if there exists  a function $\lambda :M^m\to (0, \infty)$ so that the conformal immersion $\phi: (M^{m}, {\bar g}=\lambda^{-2}g) \to (N^{n},h)$  is a biharmonic map. When this happens, we also say that the submanifold $\phi: (M^{m}, g) \to (N^{n},h)$ {\em can be biharmonically conformally immersed into its ambient space}.

We know  that a surface $\phi: (M^{2}, g) \to (N^{n},h)$ is minimal if and only if its defining isometric immersion $\phi$ a harmonic. In this case, it is well known that for any function $\lambda :M^2\to \mathbb{R}^{+}$, the conformal immersion $\phi: (M^{2}, {\bar g}=\lambda^{-2}g) \to (N^{n},h)$  is again harmonic and hence trivially biharmonic. So every minimal surface admits a trivial biharmonic conformal immersion into its ambient space. \\
\indent An example of a surface that admits a proper biharmonic conformal immersion into its ambient space is the circular cylinder 
$\phi:\mathbb{R}^2\to (\mathbb{R}^3, \delta_{ij})$,
$\phi(x, y)=(R\cos\,\frac{x}{R}, R\sin\,\frac{x}{R}, y)$ in Euclidean $3$-space. As we have seen in Proposition \ref{EZHU} that  there exists $\lambda: \mathbb{R}^2\to \mathbb{R}^+$ with $\lambda^{-2}(x, y)=e^{y/R}$ so that the conformal immersion $\phi:(\mathbb{R}^2, \lambda^{-2}({\rm d}x^2+{\rm d}y^2))\to (\mathbb{R}^3, \delta_{ij})$ is a proper biharmonic map.\\
\indent The following theorem shows that, locally,  a circular cylinder is the only one among surfaces of constant mean curvature without umbilical point that admit a proper biharmonic conformal immersion into $\mathbb{R}^3$.
\begin{theorem}\cite{Ou15}
If a surface in $\mathbb{R}^3$ with  non-zero constant mean curvature and no umbilical point  can be biharmonically conformally immersed into $\mathbb{R}^3$,  then, locally,  it is  an open portion of a circular cylinder.
\end{theorem}
\indent By a well-known results of \cite{CMO01}, any proper biharmonbic isometric immersion of a surface into $3$-sphere is given by the the standard isometric embedding $S^2(\frac{1}{\sqrt{2}})  \to S^3$ or a part of it. An interesting question to ask is whether there is a surface in Euclidean $3$-sphere that admits a proper biharmonic conformal immersion into $S^3$. The following theorem gives a partial answer to the question.
\begin{theorem}\cite{Ou17}
If $\phi: (M^2, g=\phi^*h) \to (S^3,h)$ is a compact surface in Euclidean  $3$-sphere with  $|A|^2\ge 2$, and suppose that  there is a function $\lambda$ so that the conformal immersion $\phi:(M^2, {\bar g}=\lambda^{-2} g)\to (S^3, h)$ is biharmonic. Then,  either $M^2$ is minimal, or $\lambda$ is constant, $M^2= S^2(\frac{1}{2})$, and  $\phi: M^2\to S^3$ is the standard embedding $S^2(\frac{1}{\sqrt{2}})  \to S^3$.
\end{theorem}
\indent Finally, we close this review with the following preliminary result about biharmonic conformal immersions of surfaces into a Euclidean space obtained from a consideration similar to Weierstrass representation of harmonic conformal immersions.
\begin{theorem}\label{Weie}\cite{Ou09}
A map $\phi : (M^2, g=\lambda^2({\rm d}u^2+{\rm d}v^2)) \to \mathbb{R}^n$ into Euclidean space is a  biharmonic conformal immersion if and only if
the section $ \phi=\frac{\partial \varphi}{\partial z}=\frac{1}{2}(
\varphi_{u}-i \varphi_{v})=\phi^\alpha(z)\frac{\partial}{\partial
y^\alpha}$ solves the system
\begin{eqnarray}\label{Wrep3}
&&\sum _{\alpha=1}^{n}(\phi^\alpha)^2=0,\;\;\; \sum
_{\alpha=1}^{n}|\phi^\alpha|^2\ne 0,\\
&&\frac{\partial}{\partial {\bar z}}\frac{\partial}{\partial
z}\left(\lambda^{-2}\frac{\partial \phi^\sigma}{\partial {\bar
z}}\right)=0,\;\;\;\; \sigma=1, 2,\ldots, n.
\end{eqnarray}
\end{theorem}
\begin{ack}
The author would like to thank Cezar Oniciuc and the referee for some comments and suggestions that help to improve the original manuscript. 
\end{ack}

\end{document}